\documentclass[a4paper,12pt,mingoth]{amsart}
\usepackage{amsmath,amssymb,amsthm,ascmac,fancybox}   
\usepackage{bm}              
\usepackage{tabularx}       
\usepackage[dvipdfm]{graphicx}
\allowdisplaybreaks[3]
\newcommand{\sa}   [1]{\mathcal{#1}}

\newcommand{\norm} [2]{\left\| #1 \right\|_{#2}}             
\newcommand{\ab}  [1]{\left| #1\right|}

\newcommand{\qw}   [0]{\par{\noindent}}

\newcommand{\vp}{\varphi}

\newcommand{\N}{\mathbb{N}}

\newcommand{\R}{\mathbb{R}}

\newcommand{\alm}{\mathrm{a.e.}\,}

\newcommand{\kko}  [1]{\left( #1\right)}
\newcommand{\kkko}  [1]{\left\{ #1\right\}}

\newcommand{\owari}  [0]{\rule{5pt}{10pt}}
\newcommand{\ds}{\displaystyle}

\title{The de la Vall\'{e}e Poussin mean and polynomial approximation for exponential weight}
\author{K.\ Itoh,  R.\ Sakai and N.\ Suzuki }	

\begin{document}
\maketitle
\begin{abstract}
We study $L^{p}$ boundedness of the de\ la\ Vall\'{e}e\ Poussin\ means $v_{n}(f)$ for exponential weight $w(x) = \exp(-Q(x))$ 
on $\R$. Our main result is 
\[ \left\|v_n(f)\frac{w}{T^{1/4}}\right\|_{L^p(\R)}\le C\|fw\|_{L^p(\R)}, \]
for every $n \in \N$ and every $1 \leq p \leq \infty$, where $T(x) = xQ'(x)/Q(x)$. As an application, we obtain
$\lim_{{n \rightarrow \infty}} \left\|(f - v_n(f))w/T^{1/4}\right\|_{L^p(\R)} = 0$ for 
$fw \in L^{p}(\R)$. 
\end{abstract}
\section{{Introduction}}
Let $\R=(-\infty,\infty)$ and $\N=\{1,2,\cdots,\}$. We consider an exponential weight
\[ w(x)=\exp(-Q(x)) \]
on $\R$, where $Q$ is an even and nonnegative function on $\R$.\ Throughout this
paper we always assume that $w$ belongs to a relevant class $\sa{F}(C^2+)$ (see
Section 2). We consider a function $T$ defined by
\[ T(x):=\frac{xQ'(x)}{Q(x)},~~x\neq 0. \leqno{(1.1)} \]
If $T$ is bounded, then $w$ is called a Freud-type weight, and otherwise, $w$ is
called an Erd{\H{o}}s-type weight.

The Mhaskar-Rakhmanov-Saff number (MRS number) $a_n$ for $w=\exp(-Q)$ is defined by the following equation:
\[ n=\frac{2}{\pi}\int_0^1\frac{a_nuQ^\prime(a_nu)}{(1-u^2)^{1/2}}du,~~n\in\N. \leqno{(1.2)} \]

Let $\{p_n\}$ be orthogonal polynomials for a weight $w$, that is, $p_n(x) = p_{n}(w,x)$ is the polynomial of degree $n$ such that
\[ \int_\R p_n(x)p_m(x)w^{2}(x)dx = \delta_{mn}. \]
For $1\le p\le\infty$, we denote by $L^p(\R)$ the usual $L^p$ space on $\R$.\ For a function $f$ with $fw\in L^p(\R)$, we set
\[ s_n(f)(x):=\sum_{k=0}^{n-1}c_k(f)p_k(x),~~\mbox{where}~~c_k(f) :=\int_\R f(t)p_k(t)w^{2}(t)dt \]
for $n\in\N$ (the partial sum of Fourier series). The de la Vall\'{e}e Poussin mean $v_n(f)$ of $f$ is defined by
\[ v_n(f)(x):=\frac{1}{n}\sum_{j=n+1}^{2n}s_j(f)(x). \leqno(1.3)\]

The main result of this paper is the following theorem. 

\vspace{2ex}

\qw{\bf Theorem 1.1.} Let $1\le p\le\infty$. We assume that $w\in\sa{F}(C^2+)$ satisfies
\[ T(a_n)\le c\left ( \frac{n}{a_n} \right )^{2/3} \leqno{(1.4)} \]
for some $c>0$. Then there exists a constant $C= C(w,p)>0$ such that when $T^{1/4}fw\in L^p(\R)$, then
\[ \|v_n(f)w\|_{L^p(\R)}\le C\|T^{1/4}fw\|_{L^p(\R)}, \leqno{(1.5)} \]
and when $fw\in L^p(\R)$, then
\[ \left\|v_n(f)\frac{w}{T^{1/4}}\right\|_{L^p(\R)}\le C\|fw\|_{L^p(\R)}. \leqno{(1.6)} \]

\vspace{2ex}

For Freud-type weights, it is known  the inequality
\[ \left\|v_n(f)w\right\|_{L^p(\R)}\le C\|fw\|_{L^p(\R)} \leqno{(1.7)} \]
(cf. \cite[Lemma 6]{F},  \cite[Theorem 3.4.2]{Ein03}). Note that when $T$ is bounded, then (1.4) holds evidently, so (1.7)
follows from (1.5) or (1.6) immediately. As for Erd\H{o}s-type weights, Lubinsky and Mthembu \cite{LMth} proved 
(1.5)  for $p=1$, (1.6) for $p= \infty$ and 
$ \left\|v_n(f)wT^{-1/4}\right\|_{L^p(\R)}\le C\|T^{{1/4}}fw\|_{L^p(\R)}$ for $1 < p < \infty$ under the assumption 
$T(x) = O(Q'(x)^{\varepsilon})$ for every $\varepsilon >0$. (They discussed on Cesaro means (= (C,1) means) mainly, but it easily ensure the result on de la Vall\'{e}e Poussin means.) 
We note that their assumption implies our condition (1.4) (see Remark 7.1 in Section 7). Hence our results (1.5) and (1.6) improve their result. For more general weights, 
see \cite{H}, \cite{LMache} and \cite {LMash}.

Using (1.6), we have the following result. There exists a constant  $C = C(w,p) >0$ such that for any 
$wf \in  L^p(\R)$, 
$$
\left\|(f - v_n(f)) \frac{w}{T^{1/4}}\right\|_{L^p(\R)} \le C \inf_{P \in {\mathcal{P}}_n} \|w(f - P)\|_{L^p(\R)},
\leqno(1.8)
$$
where ${\mathcal{P}}_n$ is the set of all polynomials of degree at most $n$.  It is known that the right hand side of (1.8) attains some $P_n \in {\mathcal{P}}_n$, however, 
it is not easy to determine such $P_n$ explicitly. (1.8) means that the de la Vall\'{e}e Poussin mean takes the place of $P_n$ in some sense, i.e., $v_n(f)$ 
is a good concrete approximation polynomial for given function $f$.

This paper is organized as follows. The definition of class ${\mathcal{F}}(C^2+)$ is
given in Section 2. Recalling some estimates for exponential weights in Section
3, we will give a proof of Theorem 1.1 for the case of $p = \infty$ in Section 4. 
Basic method of our proof is classic and well-known (\cite{F}, see also \cite{LMth}, \cite{LMache}), but we repeat it in order to make a role of $T(x)$ clear.
The complete proof of Theorem 1.1 is given in Section 5. 
The estimate (1.8) is shown in Section 6 together with another application. We discuss about the condition (1.4) in final Section 7.

 Throughout this paper $C$ will denote a
positive constant whose value is not necessary the same at each occurrence;
it may vary even within a line. When we write $C=C(a,b,\cdots)$, then $C$ is a
constant which depends on $a,b,\cdots$ only.

\section{{Definitions and notation}}
We say that an exponential weight $w=\exp(-Q)$ belongs to class $\sa{F}(C^2+)$, when $Q:\R\to[0,\infty)$ is a continuous and even function and satisfies the following conditions:\\
$~~$(a)$~~Q^\prime(x)$ is continuous in $\R$ with $Q(0)=0$.\\
$~~$(b)$~~Q^{\prime\prime}(x)$ exists and is positive in $\R\setminus\{0\}$.\\
$~~$(c)$\ds~~\lim_{x\to\infty}Q(x)=\infty.$\\
$~~$(d)$~~$The function $T$ defined in (1.1) is quasi-increasing in $(0,\infty)$(i.e., there exists $C>0$ such that $T(x) \le C T(y)$ whenever $0<x<y$), and there exists $\Lambda\in\R$ such that
\begin{align*} 
T(x)\ge \Lambda>1,~~x\in \R\setminus\{0\}.
\end{align*}
$~~$(e)$~~$There exists $C>0$ such that
\begin{align*} 
\frac{Q''(x)}{|Q'(x)|}\le C\frac{|Q'(x)|}{Q(x)},~~\alm~x\in \R
\end{align*}
and  there also exist a compact subinterval $J (\ni 0)$ of $\R$ and $C>0$ such that
\begin{align*}
\frac{Q''(x)}{|Q'(x)|}\ge C\frac{|Q'(x)|}{Q(x)},~~\alm~x\in \R\setminus J.
\end{align*}

\vspace{2ex}

Let $w(x)=\exp(-Q(x))\in\sa{F}(C^2+)$. Suppose that $Q\in C^3(\R)$ and $ \lambda >0$. If there exist $C >0$ and $K>0$ such that for $|x|\ge K$,
\[ \ab{\frac{Q^{\prime\prime\prime}(x)}{Q^{\prime\prime}(x)}}\le C \ab{\frac{Q^{\prime\prime}(x)}{Q^{\prime}(x)}} \ \ \mbox{and} \ \ 
\frac{|Q^\prime(x)|}{Q^\lambda(x)}\le C, \]
then we write $w\in\sa{F}_\lambda(C^3+)$. Note that if $w$ is a Freud-type, the last inequality holds with $\lambda=1$.

A typical example of Freud-type weight is $w(x) =\exp(-|x|^\alpha)$ with $\alpha>1$.  Note that  the Hermite polynomials are the orthogonal polynomials for the weight $\exp(-|x|^2)$. 
Let
$u\ge 0,\ \alpha>0,\ \alpha+u>1,\ \ell\in\N$, and we set 
\[ Q(x):=|x|^u(\exp_\ell(|x|^\alpha)-\exp_\ell (0)), \]
where $\exp_\ell(x)=\exp(\exp(\exp(\cdots(\exp x))))$ ($\ell$-times). Then $w(x) =\exp(-Q(x))$
is an Erd{\H{o}}s-type weight, which belongs to $\sa{F}_\lambda(C^3+)$ with $\lambda >1$ (see \cite{Ein01}).

We recall some notation which we use later (cf. \cite{Ein02}). By definition, the partial sum of Fourier series is given by
\[ s_m(f)(x)=\int_\R K_m(x,t)f(t)w^{2}(t)dt, \leqno(2.1)\]
where
\[ K_m(x,t) = \sum_{k=0}^{m-1}p_k(x)p_k(t). \leqno(2.2) \]
It is known that by the Christoffel-Darboux formula,
\[ K_m(x,t) = \frac{\gamma_{m-1}}{\gamma_m}\left ( \frac{p_m(x)p_{m-1}(t)-p_m(t)p_{m-1}(x)}{x-t} \right ) , \leqno(2.3)\]
where  $\gamma_n$ is the leading coefficient of $p_n$, i.e., $p_n(x) = \gamma_n x^n + \cdots $.

The Christoffel function $\lambda_{m}(x) = \lambda_m(w,x)$ is defined by
\[ \lambda_m(x):=\frac{1}{K_m(x,x)}=\kko{\sum_{k=0}^{m-1}p_k^{2}(x)}^{-1}, \leqno(2.4)\]
and $L^p$ type Christoffel function $\lambda_{m,p}(x) = \lambda_{m,p}(w,x)$ is defined by
\[ \lambda_{m,p}(x):=\inf_{P\in\sa{P}_m} \frac{1}{|P(x)|^{p}} \int_\R |P(t)w(t)|^pdt. \]
If $p=2$, we know 
\[ \lambda_{m,2}(x)=\lambda_m(x). \]

The MRS numbers $\{a_n\}$ for weight $w$ is monotone increasing, and we see easily 
\begin{align*} 
\lim_{n \to \infty} a_n = \infty ~~\textrm{and}~~ \lim_{n \to \infty} \frac{a_n}{n} = 0.
\end{align*}
Moreover, for $n\in \N$, there exists $C >0$ independent of $n$, such that
\[ \frac{a_n}{C} \le \frac{\gamma_{n-1}}{\gamma_n} \le C a_n. 
\leqno(2.5)\]
(see \cite[Lemma 13.9]{Ein02}). We also use the following functions:
$$
\Phi_n(x) :=
\begin{cases}
       \displaystyle \frac{1-\frac{|x|}{a_{2n}}}{\sqrt{1-\frac{|x|}{a_n}+\delta_n}},& |x|\le a_n, 
       \\[5ex]
         \Phi_n(a_n),& a_n<|x|, 
\end{cases}
\leqno(2.6)
$$
where $\delta_n :=\{nT(a_n)\}^{-2/3}$ 
and 
$$
\varphi_n(x):=\frac{a_n}{n}\Phi_n(x).
\leqno(2.7)
$$

\section{Lemmas}
We recall  some basic estimates. 
\vspace{2ex}

\qw{\bf Lemma 3.1.} (infinite-finite range inequality) {\cite[Theorem 1.9 (a)]{Ein02}} Let $w\in\sa{F}(C^2+)$, $0< p\le\infty$ and  $P\in\sa{P}_n\ (n\ge 1)$. Then
\[ \|Pw\|_{L^p(\R)} \le 2\|Pw\|_{L^p(|x|<a_n)}.  \leqno(3.1)\]

\vspace{2ex}

\qw{\bf Lemma 3.2.} {\cite[Lemma 3.4]{Ein05}} Let $w\in\sa{F}(C^2+)$. There exists a constant $C=C(w)>0$ such that 
\[ \frac{a_n}{n}\frac{1}{\sqrt{T(x)}\vp_n(x)}\le C, \leqno(3.2)\]
and hence 
\[ \frac{1}{\sqrt{T(x)}}\le C\Phi_n(x) \leqno(3.3)\]
holds for every $n\in\N$ and $x>0$.

\vspace{2ex}

\qw{\bf Lemma 3.3.} {\cite[Theorem 9.3]{Ein02}} Let $w\in\sa{F}(C^2+)$, and $0< p\le\infty$. Then there exists a constant $C=C(w,p)>0$ such that for every $n \in \N$ and $x \in \R$, 
\[ \lambda_{n,p}(x) \ge C\vp_n(x)w^{p}(x)\leqno(3.4) \]
and when $|x|\le a_n$, we have 
\[ \frac{1}{C}\vp_n(x)w^{p}(x) \le \lambda_{n,p}(x) \le C\vp_n(x)w^{p}(x). \leqno(3.5)\]

\vspace{2ex}

\qw{\bf Lemma 3.4.} Let $w\in\sa{F}(C^2+)$. Then there exist $0<c<1$ and $C=C(w)>0$ such that 
\[ \frac{1}{C}T(x) \le T\kko{x\pm\frac{c}{T(x)}}\le CT(x) \leqno(3.6)\]
for every $x \in \R$.\\
$[$Proof$]$ It is known that there exist $0<c_0<1$ and $C=C(w)>0$ such that
\[ \frac{1}{C}T(x) \le T\kko{x\left[1\pm\frac{c_0}{T(x)}\right]}\le CT(x) \leqno(3.7)\]
for every $x \in \R$ ({\cite[Theorem 3.2 (e)]{Ein02}}). If $|x|\le 1$, then
\[ \ab{x\pm\frac{c_0}{T(x)}}\le |x|+c_0\le 2, \]
and hence (3.6) holds evidently. If $|x|>1$, then
\[ |x|\kko{1+\frac{c_0}{T(x)}}\ge\ab{x\pm\frac{c_0}{T(x)}}\ge|x|\kko{1-\frac{c_0}{T(x)}}. \]
Since $T$ is quasi-increasing, (3.7) shows (3.6).\hfill{$\owari$}

\vspace{2ex}

\qw{\bf Lemma 3.5.} {\cite[Theorem 3.2]{Ein05}} If $w\in\sa{F}(C^2+)$ is an Erd{\H{o}}s-type weight, then for any $\eta>0$, there exists a constant $C=C(w,\eta)>0$ such that
\[ a_n\le Cn^{\eta}. \leqno(3.8)\] 
 for every $n \in \N$.

\vspace{2ex}

\qw{\bf Lemma 3.6.} {\cite[Lemma 3.5 (a), (b) and Lemma 3.11 (3.52)]{Ein02}} Let $w\in\sa{F}(C^2+)$. There exists a constant $C=C(w)>0$ such that
\[ a_{2n}\le Ca_n, \leqno(3.9) \]
\[ T(a_{2n}) \le CT(a_n)\leqno(3.10) \]
and
\[ \frac{1}{C}\frac{a_n}{T(a_n)} \le a_{2n}-a_n. \leqno(3.11) \]

\section{Proof of Theorem 1.1 for $p=\infty$}

We  begin with the following proposition. 

\vspace{2ex}

\qw{\bf Proposition 4.1.} Let $w\in\sa{F}(C^2+)$. Then there exists a constant $C=C(w)>0$ such that
\[ \frac{w^{2}(x)}{\sqrt{T(x)}}\sum_{k=0}^{n-1}p_k^{2}(x) \le C\frac{n}{a_n}\leqno(4.1)\]
for every $x \in \R$ and every $n \in \N$. 

\vspace{2ex}

$[$Proof$]$ From Lemmas 3.2 and 3.3 for $p=2$, the proposition follows immediately.{\hfill{$\owari$}}\\

The following estimate is known essentially (see \cite[Theorem 4]{LMache}, \cite[Theorem 1.2]{LMash} and \cite[Lemma 1]{H}), but we give a proof for the sake of completeness. 

\vspace{2ex}

\qw{\bf Theorem 4.2.} Let $w\in\sa{F}(C^2+)$. Then there exists a constant $C=C(w)>0$ such that 
\[ \norm{v_n(f)w\Phi_{2n}^{1/2}}{L^\infty(\R)} \le C\norm{fw}{L^\infty(\R)} \leqno(4.2) \]
for any $wf\in L^\infty(\R)$ and $n \in \N$. 

\vspace{2ex}

$[$Proof$]$ Let $x\in\R$ and $n, m \in\N$ with $n-1\le m \le 2n $. We set
\begin{align*}
g(t) &:= f(t)\chi_{(x-\frac{a_n}{n},x+\frac{a_n}{n})}(t), \\
h(t) &:= f(t)-g(t),
\end{align*}
where $\chi_{(a,b)}$ is the characteristic  function of a set  $(a,b)$. Then 
\[ v_n(f)(x)=v_n(g)(x)+v_n(h)(x). \leqno(4.3)\]
We first consider an estimate of $v_n(g)$. 
By the Schwarz inequality 
\begin{align*}
|s_m(g)(x)| &= \ab{\int_{x-\frac{a_n}{n}}^{x+\frac{a_n}{n}}f(t)K_m(x,t)w^{2}(t)dt} \\
&\le \norm{fw}{L^\infty(\R)}\int_{x-\frac{a_n}{n}}^{x+\frac{a_n}{n}}|K_m(x,t)|w(t)dt \\
&\le \norm{fw}{L^\infty(\R)}\kko{\frac{2a_n}{n}}^{1/2}\kko{\int_\R K_m^{2}(x,t)w^{2}(t)dt}^{1/2},
\end{align*}
and since 
\begin{align*}
\int_\R K_m^{2}(x,t)w^{2}(t)dt & = \sum_{i=0}^{m-1}p_i^{2}(x) \le \sum_{i=0}^{2n-1}p_i^{2}(x) = \frac{1}{\lambda_{2n}(x)},
\end{align*}
(2.7) and (3.4) give us 
\begin{align*}
|w(x)s_m(g)(x)\Phi_{2n}^{1/2}(x)|\le\norm{fw}{L^\infty(\R)}\kko{\frac{2\vp_{2n}(x)}{\lambda_{2n}(x)}w^{2}(x)}^{1/2}\le C\|fw\|_{L^\infty(\R)},
\end{align*}
so that 
\[ |v_n(g)(x)w(x)\Phi_{2n}^{1/2}(x)| \le C\|fw\|_{L^\infty(\R)}. \leqno(4.4)\]

Next, for an estimate of $v_n(h)(x)$, we set 
\[ H(t):=\frac{h(t)}{x-t}, \]
and denote by $\{b_j(H)\}$ the Fourier coefficients of $H$, i.e., 
\[ b_j(H):=\int_\R H(t)p_j(t)w^{2}(t)dt. \]
Using the Christoffel-Darboux formula (2.3), we have 
\begin{align*}
|s_m(h)(x)|= &\ab{\int_\R\frac{\gamma_{m-1}}{\gamma_m}\left ( \frac{p_m(x)p_{m-1}(t)-p_m(t)p_{m-1}(x)}{x-t} \right ) h(t)w^{2}(t)dt} \\
= &\frac{\gamma_{m-1}}{\gamma_{m}}\ab{p_m(x)\int_\R p_{m-1}(t)H(t)w^{2}(t)dt-p_{m-1}(x)\int_\R p_{m}(t)H(t)w^{2}(t)dt}\\
\le & \frac{\gamma_{m-1}}{\gamma_{m}}( |p_m(x)||b_{m-1}(H)|+|p_{m-1}(x)||b_m(H)|).
\end{align*}
By (2.5) and the Schwarz inequality, we have 
\begin{align*}
|v_n(h)(x)| &\le \frac{1}{n}\max_{m\le 2n}\kkko{\frac{\gamma_{m-1}}{\gamma_{m}}}\times\frac{2}{\lambda_{2n}^{1/2}(x)}\kko{\sum_{m=0}^{2n}b_m^{2}(H)}^{1/2} \\
&\le C\frac{a_n}{n}\frac{1}{\lambda_{2n}^{1/2}(x)}\kko{\sum_{m=0}^{2n}b_m^{2}(H)}^{1/2}.
\end{align*}
 Since 
\begin{align*}
\sum_{m=0}^{2n}b_m^2(H) &\le \sum_{m=0}^{\infty}b_m^2(H) \le \int_\R(H(t)w(t))^2dt \\
&= \int_{|t-x|\ge\frac{a_n}{n}}\kko{\frac{f(t)w(t)}{t-x}}^2dt \le \|fw\|^2_{L^\infty(\R)}\frac{2n}{a_n},
\end{align*}
(2.7) and (3.4) imply 
\begin{align*}
|v_n(h)(x)w(x)| &\le C\frac{a_n}{n}\kko{\frac{w^{2}(x)}{\lambda_{2n}(x)}}^{1/2}\norm{fw}{L^\infty(\R)}\sqrt{\frac{n}{a_n}} \\
&\le C\norm{fw}{L^\infty(\R)}\frac{1}{\Phi_{2n}^{1/2}(x)}\kko{\frac{w^{2}(x)\vp_{2n}(x)}{\lambda_{2n}(x)}}^{1/2} \\
&\le C\norm{fw}{L^\infty(\R)}\frac{1}{\Phi_{2n}^{1/2}(x)}.
\end{align*}
This together with (4.4) gives us 
\begin{align*}
|v_n(f)(x)w(x)\Phi_{2n}^{1/2}(x)| \le C\norm{fw}{L^\infty(\R)},
\end{align*}
which completes the proof of (4.2). {\hfill{$\owari$}} 

\vspace{2ex}

The following corollary contains (1.6) for $p = \infty$. 

\vspace{2ex}

\qw{\bf Corollary 4.3.} Let $w\in\sa{F}(C^2+)$. Then there exists a constant $C=C(w)>0$ such that
\[ \norm{v_n(f)\frac{w}{T^{1/4}}}{L^\infty(\R)} \le C\norm{fw}{L^\infty(\R)} \]
 for any $wf\in L^\infty(\R)$.
 
 \vspace{2ex}
 
$[$Proof$]$ This follows from Theorem 4.2 and Lemma 3.2  immediately.{\hfill{$\owari$}}\\

\vspace{2ex}

We give a proof of (1.5) for $p = \infty$.

\vspace{2ex}

\qw{\bf Theorem 4.4.} Let $w\in\sa{F}(C^2+)$ such that  $T(a_n)\le c(n/a_n)^{2/3}$ with some $c >0$. Then there exists a constant $C = C(w) >0$ such that if $T^{1/4}wf\in L^\infty(\R)$, we have
\[ \norm{v_n(f)w}{L^\infty(\R)} \le C \norm{T^{1/4}fw}{L^\infty(\R)}. \leqno(4.5) \]

\vspace{2ex}

$[$Proof$]$ Let $x\in\R$ and $n, m \in\N$ such that $n-1\le m \le 2n$. As in the proof of Theorem 4.2,  we set $
g(t) := f(t)\chi_{(x-\frac{a_n}{n},x+\frac{a_n}{n})}(t), \ 
h(t) := f(t)-g(t).$
Then $ v_n(f)(x)=v_n(g)(x)+v_n(h)(x)$. 
Since $v_n(g), \ v_n(h) \in {\mathcal{P}}_{2n}$, by Lemma 3.1, we may suppose that $|x|\le a_{2n}$.

We first prove that for any $n\in\N$,
\[ T(x)\le CT(t) \leqno{(4.6)} \]
whenever $|x|\le a_{2n}$ and $|x-t|\le a_n/n$. In fact, by (1.4) and (3.11), 
there exists $N_0\in\N$ such that for every $n\ge N_0$,  
\[ |t|\le|x|+\frac{a_n}{n}\le a_{2n}+\frac{1}{C} \frac{a_n}{T(a_n)}\le a_{4n}, \]
which shows $ T(t) \leq CT(a_{4n}) \le CT(a_n) \leq C (n/a_n)^{2/3} \leq c_0n/a_n$ by (3.10), 
where $c_0$ is the constant in Lemma 3.4. Hence for any $n\ge N_0$,
\[ |x|\le \frac{a_n}{n}+|t| \le \frac{c_0}{T(t)}+|t|, \] 
so that Lemma 3.4 gives us 
\begin{align*}
T(x) = T(|x|) \le CT\kko{\frac{a_n}{n}+|t|} 
\le CT\kko{\frac{c_0}{T(t)}+|t|}\le CT(|t|)=CT(t).
\end{align*}
 When $n\le N_0$, then $1\le T(x)\le C$ for $|x|\le a_{2N_0}$, so that (4.6) holds immediately.

Now we estimate $v_n(g)(x)$. By the Schwarz inequality, (4.6) and Proposition 4.1,
\begin{align*}
|s_m(g)(x)w(x)| 
&\le C\kko{ \sum_{k=0}^{m-1}\frac{w^{2}(x)}{T^{1/2}(x)}p_k^{2}(x)}^{1/2}\kko{\int_{x-\frac{a_n}{n}}^{x+\frac{a_n}{n}}|T^{1/4}(t)w(t)f(t)|^2dt}^{1/2} \\
&\le C\sqrt{\frac{m}{a_m}}\norm{T^{1/4}wf}{L^\infty(\R)}\kko{\int_{x-\frac{a_n}{n}}^{x+\frac{a_n}{n}}dt}^{1/2} \\
&\le C\norm{T^{1/4}wf}{L^\infty(\R)},
\end{align*}
which implies 
\[ |v_n(g)(x)w(x)|\le C\norm{T^{1/4}fw}{L^\infty(\R)}.
\leqno(4.7) \] 

Next, we estimate $v_n(h)(x)$.\ Since 
\[ v_n(h)(x)=\frac{1}{n}\sum_{m=n+1}^{2n}\frac{\gamma_{m-1}}{\gamma_{m}}(b_{m-1}(H)p_m(x)-b_m(H)p_{m-1}(x))\]
as before, where $H(t) = h(t)/(x-t)$, Proposition 4.1 gives us 
\begin{align*}
|w(x)v_n(h)(x)| & \le \frac{2a_n}{n}\kko{\sum^{2n}_{{m=0}} w^{2}(x) \sum^{2n}_{m=n} p_{m}^{2}(x)}^{1/2}\kko{\sum_{m=0}^{2n}b_m^{2}(H)}^{1/2} \\
& \le C\sqrt{\frac{a_n}{n}} \left ( T^{1/2}(x)\int_{|x-t|>\frac{a_n}{n}}\ab{\frac{f(t)w(t)}{x-t}}^2dt \right )^{1/2}.
\end{align*}
By Lemma 3.4, we can take $c>0$ small enough such that $T(x)\ge C T(x\pm c/T(x))$ for $C>0$.\ Then we have
\begin{align*}
& T^{1/2}(x)\int_{\frac{c}{T(x)}>|x-t|>\frac{a_n}{n}}\frac{|f(t)w(t)|^2}{(x-t)^2}dt \\
\ & \le  CT^{1/2}(x)\int_{\frac{c}{T(x)}>|x-t|>\frac{a_n}{n}}\frac{1}{T^{1/2}(|x-\frac{c}{T(x)}|)}\frac{|f(t)T^{1/4}(t)w(t)|^2}{(x-t)^2}dt\\
\ & \le C\|T^{1/4}fw\|_{L^\infty(\R)}^2\int_{\frac{c}{T(x)}>|x-t|>\frac{a_n}{n}}\frac{1}{(x-t)^2}dt \\
\ & \le  C\frac{n}{a_n}\|T^{1/4}fw\|_{L^\infty(\R)}^2.
\end{align*}
Also since $T(a_n)\le C (n/a_{n})^{2/3}$ and $|x|\le a_{2n}$, we have
\begin{align*}
& T^{1/2}(x)\int_{\frac{c}{T(x)}\le |x-t|}\frac{|f(t)w(t)|^2}{(x-t)^2}dt \\
\ \ &  \le\|fw\|_{L^\infty(\R)}^2T^{1/2}(x)\int_{\frac{c}{T(x)}\le |x-t|}\frac{1}{(x-t)^2}dt \\
\ \ & \le C\|fw\|_{L^\infty(\R)}^2T^{3/2}(x)\le C \frac{n}{a_{n}}\|fw\|_{L^\infty(\R)}^2.
\end{align*}
These estimates show 
\begin{align*}
\kko{T^{1/2}(x) \int_{|x-t|>\frac{a_n}{n}}\ab{\frac{f(t)w(t)}{x-t}}^2dt}^{1/2}\le C\sqrt{\frac{n}{a_n}}\|T^{1/4}fw\|_{L^\infty(\R)},
\end{align*}
and hence 
\begin{align*}
\|v_n(h)w\|_{L^\infty(\R)} \le C\|T^{1/4}fw\|_{L^\infty(\R)}
\end{align*}
follows.  This together with (4.7) implies (4.4). {\hfill{$\owari$}}

\section{Proof of Theorem 1.1 for $1\le p\le\infty$}

In this section we complete the proof of Theorem 1.1.

\vspace{2ex}

$[$Proof of (1.5) $]$ The $L^\infty$-norm case is Theorem 4.4. We prove the $L^1$-norm case. By the duality of $L^1$-norm,
\[ \|v_n(f)w\|_{L^1(\R)}=\sup_{\|gw\|_{L^\infty(\R)}\le 1}\int_\R v_n(f)(x)q(x)w(x)^2dx. \]
Since  $K_m(x,t)=K_m(t,x)$, we see
\[ \int_\R s_m(f)(x)g(x)w^2(x)dx = \int_\R f(x)s_m(g)(x)w^2(x)dx,\]
and hence 
\[ \int_\R v_n(f)(x)g(x)w^2(x)dx = \int_\R f(x)v_n(g)(x)w^2(x)dx. \]
Therefore, using Corollary 4.3, we have 
\begin{align*}
\|v_n(f)w\|_{L^1(\R)} &= \sup_{\|gw\|_{L^\infty(\R)}\le 1}\int_\R v_n(f)(x)g(x) w^2(x)dx \\
&= \sup_{\|gw\|_{L^\infty(\R)}\le 1}\int_\R f(x)v_n(g)(x)w^{2}(x)dx \\
&\le \sup_{\|gw\|_{L^\infty(\R)}\le 1}\norm{T^{1/4}fw}{L^1(\R)}\norm{\frac{1}{T^{1/4}}v_n(g)w}{L^\infty(\R)} \\
&\le \sup_{\|gw\|_{L^\infty(\R)}\le 1}\norm{T^{1/4}fw}{L^1(\R)}\norm{gw}{L^\infty(\R)} \\
&\le \norm{T^{1/4}fw}{L^1(\R)}.
\end{align*}
Since the operator  norms 
\[ F \mapsto\ wv_n\kko{F\frac{1}{wT^{1/4}}}, \]
 for $p=1$ and $p=\infty$ are bounded,  the Riesz-Thorin interpolation theorem gives us 
\[ \sup_{\|F\|_{L^p(\R)\le 1}}\norm{wv_n\kko{F\frac{1}{wT^{1/4}}}}{L^p(\R)} < \infty \]
for every $p$ with $1 \leq p \leq \infty$. This implies   (1.5). \hfill{$\owari$}\\

$[$Proof of (1.6) $]$  The $L^\infty$-norm case is Corollary 4.3. Now, we show $L^1$-norm case. Similarly as above, 
\begin{align*}
\norm{v_n(f)\frac{w}{T^{1/4}}}{L^1(\R)} &= \sup_{\|wg\|_{L^\infty(\R)=1}}\int_\R v_n(f)(x)\frac{g(x)}{T(x)^{1/4}}w^2(x)dx \\
&= \sup_{\|wg\|_{L^\infty(\R)=1}}\int_\R f(x)v_n\kko{\frac{g}{T^{1/4}}}(x)w^2(x)dx \\
&\le  \sup_{\|wg\|_{L^\infty(\R)=1}} \norm{wv_n\kko{\frac{g}{T^{1/4}}}}{L^\infty(\R)}\norm{fw}{L^1(\R)}.
\end{align*}
Since 
$$
\norm{wv_n\kko{\frac{g}{T^{1/4}}}}{L^\infty(\R)} \le C\norm{T^{1/4}\frac{g}{T^{1/4}}w}{L^\infty(\R)} = C\|gw\|_{L^\infty(\R)}=C
$$
by (4.5), we see
\begin{align*}
\norm{v_n(f)\frac{w}{T^{1/4}}}{L^1(\R)} \le C\|fw\|_{L^1(\R)}. 
\end{align*}
Hence by the Riesz-Thorin interpolation theorem for  the operator
\[ F\ \mapsto\ \frac{w}{T^{1/4}}v_n(w^{-1}F), \]
we have 
\[ \sup_{\|F\|_{L^p(\R)\le 1}}\left\|\frac{w}{T^{1/4}}v_n(w^{-1}F)\right\|_{L^{p}(\R)} < \infty \]
for every $1 \leq p \leq \infty$. This implies (1.6). \hfill{$\owari$}

\section{Corollaries}

We begin with the following corollary.

\vspace{2ex}

\qw{\bf Corollary 6.1.}  Let $1\le p \le \infty$  and $w\in\sa{F}(C^2+)$ satisfy (1.4). Then there exists a constant $C=C(w,p)>0$ such that 
for every $wf\in L^p(\R)$ and $n \in \N$,  
\[  \norm{v_n(f)w}{L^p(\R)} \le C T^{1/4}(a_n) \norm{fw}{L^p(\R)}. 
\leqno(6.1) \]

\vspace{2ex}

$[$Proof$]$ By Lemma 3.1 and (1.6) in Theorem 1.1, 
\begin{align*}
\|v_n(f)w \|_{L^p(\R)} &\le C\|v_n(f)w\|_{L^p([-a_{2n},a_{2n}])} \\
&\le CT^{1/4}(a_n)\norm{v_n(f)\frac{w}{T^{1/4}}}{L^p(\R)} \\
&\le CT^{1/4}(a_n)\|fw\|_{L^p(\R)}. \tag*{$\owari$}
\end{align*}

To discuss  polynomial approximations, we define the degree of weighted polynomial approximation for $wf \in L^p(\R)$ by 
\[ E_{p,n}(w,f):=\inf_{P\in \sa{P}_n}\|w(f-P)\|_{L^p(\R)}. 
\leqno(6.2) \]  
We quote two results from our previous papers. 

\vspace{2ex}

\qw{\bf Proposition 6.2.} (\cite[Theorem 1]{Ein04} and \cite[Theorem 6.1]{Ein05})  Let $1\le p \le \infty$  and $w\in\sa{F}(C^{2}+)$. Then there exists a 
constant $C = C(w,p) >0$ such that for every $n \in N$, if $f$ is absolutely continuous and  $f'w \in L^{p}(\R)$, then 
$$E_{p,n}(w,f) \leq C \frac{a_{n}}{n}  \|f'w\|_{L^{p}(\R)} \leqno(6.3)$$
and if $P \in {\mathcal{P}}_{n}$, then 
\[ \norm{\frac{1}{\sqrt{T}}P^\prime w}{L^p(\R)}\le C\frac{n}{a_n}\|Pw\|_{L^p(\R)} \leqno(6.4) \]
holds true.

\vspace{2ex}

\qw{\bf Corollary 6.3.} Let $1\le p \le \infty$  and $w\in\sa{F}(C^2+)$ satisfy (1.4). 
Then there exists a constant $C=C(w,p)>0$ such that  for every $n \in \N$  and every $wf\in L^p(\R)$, 
\[ \norm{(f-v_n(f))\frac{w}{T^{1/4}}}{L^p(\R)}\le CE_{p,n}(w,f), \leqno(6.5) \]
\[ \norm{(f-v_n(f))w}{L^p(\R)}\le CT^{1/4}(a_n)E_{p,n}(w,f), \leqno(6.6) \]
and when $T^{1/4}wf\in L^p(\R)$, 
\[ \norm{(f-v_n(f))w}{L^p(\R)}\le CE_{p,n}(T^{1/4}w,f). \leqno(6.7) \]
Moreover if $f$ is absolutely continuous and $f'w \in L^{p}(\R)$, then 
\[ \norm{(f-v_n(f))\frac{w}{T^{1/4}}}{L^p(\R)}\le C \frac{a_{n}}{n} \|f'w\|_{L^{p}(\R)}. \leqno(6.8) \]

\vspace{2ex}

$[$Proof$]$ Since $v_n(P) = P$ for every $P \in {\mathcal{P}}_n$, we have $v_n(f)=P+v_n(f-P)$, and hence (1.6) gives us 
\begin{align*}
\norm{(f-v_n(f))\frac{w}{T^{1/4}}}{L^p(\R)} &\le \norm{(f-P)\frac{w}{T^{1/4}}}{L^p(\R)}+\norm{v_n(f-P)\frac{w}{T^{1/4}}}{L^p(\R)} \\
&\le \norm{(f-P)\frac{w}{T^{1/4}}}{L^p(\R)}+\norm{(f-P)w}{L^p(\R)}\\
&\le C\norm{(f-P)w}{L^p(\R)},
\end{align*}
which shows (6.5). Using (6.1) and (1.5), we also obtain (6.6) and (6.7), respectively. Favard type inequality 
(6.8) follows from (6.3) and (6.5).  {\hfill{$\owari$}}

\vspace{2ex}

Combining  (1.5) with (6.4), we obtain the following estimate of the derivative of $v_{n}(f)$:

\vspace{2ex}

\qw{\bf Corollary 6.4.} Let $1\le p \le \infty$  and $w\in\sa{F}(C^2+)$ satisfy (1.4).
Then there exists a constant $C=C(w,p)>0$ such that  for every $n \in \N$ 
$$
\norm{\frac{1}{\sqrt{T}}v_n^\prime(f) w}{L^p(\R)} \le C\frac{n}{a_{n}}\|T^{1/4}fw\|_{L^p(\R)}.
\leqno(6.9)$$

\section{Remarks}

In this final section we make three remarks. First one is concerning to the condition (1.4) in Theorem 1.1.

\vspace{2ex}

\qw{\bf Remark 7.1.} Let $w =\exp(-Q) \in \sa{F}(C^{2} +)$. 

(1) Let $0 < \lambda < 2$. If $w$ satisfies 
$$
|Q'(x)|/Q^{\lambda}(x) \leq C \ \ (|x| \geq 1) \leqno(7.1)
$$
for some $C >0$.  Then the condition (1.4) holds true. In particular, all $w \in \sa{F}_{\lambda}(C^{3})$ satisfy (1.4). 

(2) If $T(x) = O((Q'(x))^{\varepsilon})$ for some $0 < \varepsilon < 1/2$, then (7.1) holds for some $0 < \lambda < 2$. In particular, all the weights discussed in \cite{LMth} satisfy (1.4).  

In fact, if $w$ is Freud-type, then (1.4) holds clearly, so that we may assume that $w$ is Erd\H{o}s-type. In \cite[Lemma 3.2]{Ein05}, we showed that when $w$ satisfies (7.1), then  $T(a_{n}) \leq C n^{2/3 - \delta}$ holds for some $0 < \delta < 2/3$. Hence this and (3.8) imply (1.4). Since 
$|xQ'(x)|/Q(x) = T(x) \leq C |Q'(x)|^{\varepsilon}$, $\lambda := 1/(1-\varepsilon) <2$ satisfies (7.1).

\vspace{2ex}
Next  we remark on the degree of weighted polynomial approximation $E_{p,n}(w,f)$. 

\vspace{2ex}

\qw{\bf Remark 7.2.} Let $1 \leq p \leq \infty$.  It is known that if $w \in \sa{F}(C^{2} +)$, then $\lim_{n \rightarrow \infty} E_{p,n}(w,f) =0$ for $wf \in L^p(\R)$ (when $p=\infty$ we further assume that $\lim_{|x| \rightarrow \infty} f(x)w(x) =0)$ (e.g., \cite[Theorem 1.4]{L}). Hence (6.5) implies 
$$
\lim_{n \rightarrow \infty} \left\|(f - v_{n}(f)) \frac{w}{T^{1/4}}\right\|_{L^{p}(\R)} = 0. \leqno(7.2)
$$
This is a concrete polynomial approximation for a given function $f$. Similar argument can not apply to (6.7), because 
$T^{1/4}w$ may not belong to $\sa{F}(C^{2} +)$. To overcome this difficulty, we use a mollification of a weight (\cite[Theorem 4.1]{Ein05}): Let $w \in \sa{F}_{\lambda}(C^{3} +)$ with $0 < \lambda \leq 3/2$. Then we can construct a new weight $w^{*} \in 
\sa{F}(C^{2}+)$ which satisfies $T^{1/4}w \sim w^{*}$, $a_{n} \sim a^{*}_{n}$ and $T \sim T^{*}$, where 
$a^{*}_{n}$ and $T^{*}$ are the MRS number and  a function defined (1.1) with respect to $w^{*}$, respectively.
Hence $E_{p,n}(T^{1/4}w,f) \leq C E_{p,n}(w^{*},f)$ and (6.7) show that if $w \in \sa{F}_{\lambda}(C^{3}+)$ with $0 < \lambda \leq 3/2$ and if $T^{1/4}wf \in L^{p}(\R)$, 
 then 
 $$
\lim_{n \rightarrow \infty} \|(f - v_{n}(f)) w\|_{L^{p}(\R)} = 0 \leqno(7.3)
$$
(when $p=\infty$ we further assume that $\lim_{|x| \rightarrow \infty} T^{1/4}(x)f(x)w(x) =0)$.

\vspace{2ex}

\qw{\bf Remark 7.3.} (6.9) suggests that the following inequality would be true:
$$
\|v'_n(f) w\|_{L^p(\R)} \leq C \frac{n}{a_n} \|T^{3/4} fw\|_{L^p(\R)}.
\leqno(7.4)
$$
We will discuss this estimate elsewhere.

\mbox{}\\
{\footnotesize{
\qw Kentaro Itoh\\
Department of Mathematics\\
Meijo University\\
Tenpaku-ku, Nagoya 468-8502\\
Aichi Japan\\
{\texttt{133451501@ccalumni.meijo-u.ac.jp}}
\vspace{1ex}
\qw Ryozi Sakai\\
Department of Mathematics\\
Meijo University\\
Tenpaku-ku, Nagoya 468-8502\\
Aichi Japan\\
{\texttt{ryozi@crest.ocn.ne.jp}}
\vspace{1ex}
\qw Noriaki Suzuki\\
Department of Mathematics\\
Meijo University\\
Tenpaku-ku, Nagoya 468-8502\\
Aichi Japan\\
{\texttt{suzukin@meijo-u.ac.jp}}
}}
\end{document}